\documentclass[11pt,oneside,english]{amsart}
\usepackage[T1]{fontenc}
\usepackage[latin9]{inputenc}
\usepackage{amssymb}

\makeatletter
\numberwithin{equation}{section} 
\numberwithin{figure}{section} 
  \theoremstyle{plain}
  \newtheorem{thm}{Theorem}[section]
  \theoremstyle{definition}
  \newtheorem{defn}[thm]{Definition}
 \theoremstyle{definition}
  \newtheorem{example}[thm]{Example}
  \theoremstyle{plain}
  \newtheorem{prop}[thm]{Proposition}
  \theoremstyle{plain}
  \newtheorem{cor}[thm]{Corollary}
  \theoremstyle{remark}
  \newtheorem*{rem*}{Remark}
  \theoremstyle{remark}
  \newtheorem{rem}[thm]{Remark}

\@addtoreset{equation}{section}

\usepackage{babel}
\makeatother

\begin{document}

\title{Weyl substructures and compatible linear connections}

\author{Oana Constantinescu and Mircea Crasmareanu}

\begin{abstract}
The aim of this paper is to study from the point of view of linear
connections the data $\left(M,\mathcal{D},g,W\right),$ with $M$
a smooth $\left(n+p\right)$ di\-men\-sio\-nal real manifold, $\left(\mathcal{D},g\right)$
a \textit{$n$ }\textit{\emph{dimensional semi-Riemannian distribution}}\emph{
}on $M,$ $\mathcal{G}$ the conformal structure generated by $g$
and $W$ a Weyl substructure: a map $W:$ $\mathcal{G}\rightarrow$
$\Omega^{1}\left(M\right)$ such that $W\left(\overline{g}\right)=W\left(g\right)-du,$
$\overline{g}=e^{u}g;u\in C^{\infty}\left(M\right).$ Compatible linear
connections are introduced as a natural extension of similar notions
from Riemannian geo\-me\-try and such a connection is unique if
a symmetry condition is imposed. In the foliated case the local expression
of this unique connection is obtained. The notion of Vranceanu connection
is introduced for a pair (Weyl structure, distribution) and it is
computed for the tangent bundle of Finsler spaces, particularly Riemannian,
choo\-sing as distribution the vertical bundle of tangent bundle
projection and as 1-form the Cartan form. 
\end{abstract}



\thanks{Supported by  Grant ID 398 from the Romanian Ministery of Education. }

\date{May 4, 2009}

\maketitle

\section*{Introduction}

After Einstein's approach to gravitation \cite{Einstein}, several
others theories have been developed as part of the efforts to cure
problems arising when the gravitational field is coupled to matter
fields. Thus, as soon as Einstein presented the General Relativity,
Weyl \cite{Weyl}-\cite{w:h} proposed a new geometry in which a new
scalar field accompanies the metric field and changes the scale of
length measurements. The aim was to unify gravitation and electromagnetism,
but this theory was briefly refuted by Einstein because the non-metricity
had direct consequences over the spectral lines of elements which
has never been observed.

Let $\mathcal{G}$ be a conformal structure on the smooth manifold
$M_{n}$ i.e. an equi\-valence class of Riemannian metrics: $g\sim~\overline{g}$
if there exists a smooth function $f\in C^{\infty}\left(M\right)$
such that $\overline{g}=e^{f}g$. Denoting by $\Omega^{1}(M)$ the
$C^{\infty}\left(M\right)$-module of 1-forms on $M$, a (Riemannian)
Weyl structure is a map $W:\mathcal{G}\rightarrow\Omega^{1}\left(M\right)$
such that $W\left(\overline{g}\right)=W\left(g\right)-df$. In \cite{f:l}
and \cite{Sen} it is proved that for a Weyl manifold $\left(M,\mathcal{G},W\right)$
there exists a unique torsion-free linear connection $\nabla$ on
$M$ such that for every $g\in\mathcal{G}$ \cite{h:w}: \begin{equation}
\nabla g+W\left(g\right)\otimes g=0\label{eq:1}\end{equation}
 called \textit{Weyl connection}. The parallel transport induced by
$\nabla$ preserves the given conformal class $\mathcal{G}$. Also,
the above theory can be expressed in terms of $G$-structures with
$G$ the conformal group $CO(n)=O(n)\times\mathbb{R}^{+}$.

The literature on Weyl structures is huge and the increasing interest
in it is motivated in the last years by a new relationship with physics
and gauge theory through the notion of \textit{Weyl-Einstein manifold},
\cite{p:g}. Also, some interesting extensions of Weyl structures
inspired by generalizations of Riemannian metrics have appeared: for
Finsler metrics in \cite{a:i}, \cite{k:l} and \cite{l:k} while
for generalized Lagrange metrics in \cite{m:c}. The infinite-dimensional
case was treated in \cite{m:a}.

In this paper we propose another extension of Weyl structures and
compa-tible connections (\ref{eq:1}) namely in the semi-Riemannian
distributions framework. So, the first section is devoted to the proposed
generalization and exactly as in the semi-Riemannian geometry the
uniqueness of the compa-tible connection is obtained provided a symmetry
condition holds. Also, the compatibility condition is rewritten in
terms of quasi-connections. The second section deals with the foliated
case through the local expression in an adapted frame and a new characterization
of bundle-like metrics is obtained in terms of Weyl structures.

For Weyl structures on a manifold endowed with a distribution we introduce
the notion of Vranceanu connection following a similar tool from the
geometry of a pair (Riemannian manifold, distribution). This way,
we obtain a generalization for some notions, results and relations
of \cite{b:f}. The global expression of this connection appears in
the first section, while the local coefficients are given again for
the foliated manifolds in the second part of the paper. Let us point
out that we treat this connection considering global Weyl structures
as examples of our theory.

We devote the third section to the Vranceanu linear connection for
the tangent bundle of Finsler, particularly Riemannian spaces when
we call it \textit{Vranceanu-Cartan}, choosing as distribution the
vertical bundle of tangent bundle projection and as 1-form the Cartan
form. The motivation for this name consists in the fact that Vranceanu
in 1926 (\cite{g:v1}) and Cartan in 1928 (\cite{e:c}) are the firsts
who proposed a geometrization of non-holonomic mechanics (in the same
year, the papers \cite{j:a:s,s:y,g:v2} are devoted to the subject)
bur recently new linear connections are proposed for this framework
in \cite{b:j} and \cite{d:g}. As final problems, the flatness of
the Vranceanu-Cartan connection and the covariant derivative of the
Liouville vector fields with respect to the Vranceanu connection are
discussed.

At the end of these remarks let us point out that our work can also
be considered as proposing a generalization of the sub-Riemannian
geometry, \cite{r:m}. So, we adress a new theory, namely \emph{sub-Weyl
theory}, to which it seems to belongs also the paper \cite{z:p}.

\section{Weyl substructures and compatible connections}

For a real manifold $M$ we use the following notations: 

\begin{itemize}
\item $C^{\infty}(M)$ is the ring of smooth real functions,
\item $\chi(M)$ is  the $C^{\infty}(M)$ - module of vector fields on $M.$
\end{itemize}
\medskip{}

Let $M$ be a smooth $\left(n+p\right)$-di\-men\-sio\-nal real
manifold and $\mathcal{D}$ an $n$-di\-men\-sio\-nal distribution
on $M$. Suppose $g$ is a semi-Riemannian metric on $\mathcal{D}$,
that is, in the words of \cite[p. 23]{b:f}, $\left(\mathcal{D},g\right)$
is a \textit{semi-Riemannian distribution} on $M$. Let $\mathcal{G}=\left\{ \overline{g}=e^{u}g;u\in C^{\infty}\left(M\right)\right\} $
be the conformal structure generated by $g$.

\medskip{}

\begin{defn}
A \textit{Weyl substructure} is a map $W:\mathcal{G}\rightarrow\Omega^{1}\left(M\right)$
such that: \begin{equation}
W\left(\overline{g}\right)=W\left(g\right)-du.\label{eq:2}\end{equation}
 The data $\left(M,\mathcal{D},g,W\right)$ will be called a \textit{sub-Weyl
manifold}.
\end{defn}
\medskip{}

Let us point out that a straightforward computation gives: \[
W\left(e^{v}\overline{g}\right)=W\left(\overline{g}\right)-dv.\]
 It follows that if for some $g\in\mathcal{G}$ the 1-form $W\left(g\right)$
is closed (or exact) then for every $\overline{g}\in\mathcal{G}$
the 1-form $W\left(\overline{g}\right)$ is closed (or exact).

\medskip{}

We want a linear connection on $\mathcal{D}$ whose properties are
similar of those of Weyl connection on a Riemannian manifold. To this
end we consider a complementary distribution $\mathcal{D}^{\prime}$
to $\mathcal{D}$ in $TM$: \begin{equation}
TM=\mathcal{D}\oplus\mathcal{D}^{\prime}.\label{eq:3}\end{equation}
 Since $M$ is supposed to be paracompact there exists such a distribution.
Let $Q$ and $Q^{\prime}$ be the corresponding projectors of this
decomposition. Recall that a linear connection $\nabla$ on $\mathcal{D}$
is said to be $\mathcal{D}^{\prime}$-\textit{torsion free} if its
$\mathcal{D}^{\prime}$-torsion field vanishes i.e. \cite[p. 23]{b:f}:
\begin{equation}
\nabla_{X}QY=\nabla_{QY}QX+Q\left[X,QY\right],\quad\forall X,Y\in\mathcal{X}\left(M\right).\label{eq:4}\end{equation}

\begin{defn}
$\nabla$ is \textit{compatible} to the Weyl substructure if: \begin{equation}
\nabla_{QX}g+W\left(g\right)\left(QX\right)\cdot g=0,\quad\forall X\in\mathcal{X}\left(M\right).\label{eq:5}\end{equation}

\end{defn}
Again it results that this relation has a geometrical meaning since:
\[
\nabla_{X}\overline{g}+W\left(\overline{g}\right)\left(X\right)\overline{g}=e^{u}(\nabla_{X}g+W\left(g\right)\left(X\right)g),\quad\forall X\in\mathcal{X}\left(M\right).\]

The aim of this section is to obtain a generalization of the results
from Introduction:

\medskip{}

\begin{thm}
Given a sub-Weyl manifold with a complementary distribution $\mathcal{D}^{\prime}$\ there
exists an unique compatible linear connection $\nabla$\ on $\mathcal{D}$\ such
that $\nabla$\ is $\mathcal{D}^{\prime}$-torsion free.
\end{thm}
\medskip{}

\begin{proof}
Let us consider $\nabla$ given by: \begin{eqnarray}
2g\left(\nabla_{QX}QY,QZ\right) & = & QX\left(g\left(QY,QZ\right)\right)+QY\left(g\left(QZ,QX\right)\right)\nonumber \\
 &  & -QZ\left(g\left(QX,QY\right)\right)+g\left(Q\left[QX,QY\right],QZ\right)\nonumber \\
 &  & -g\left(Q\left[QY,QZ\right],QX\right)+g\left(Q\left[QZ,QX\right],QY\right)\nonumber \\
 &  & +W\left(g\right)\left(QX\right)g\left(QY,QZ\right)+W\left(g\right)\left(QY\right)g\left(QZ,QX\right)\nonumber \\
 &  & -W\left(g\right)\left(QZ\right)g\left(QX,QY\right)\label{eq:6}\end{eqnarray}
 respectively: \begin{equation}
\nabla_{Q^{\prime}X}QY=Q\left[Q^{\prime}X,QY\right].\label{eq:7}\end{equation}
 It is easy to verify that $\nabla$ is the unique linear connection
on $\mathcal{D}$ that satisfies the conclusion. 
\end{proof}
\medskip{}

For the particular case $p=0$ the result of \cite{f:l} and \cite{Sen}
from the Introduction is recovered.

\medskip{}

\begin{example}
Let $(M,g,W)$ be a Weyl manifold \cite{f:l} i.e. $g$ is a global
semi-Riemannian metric on $M$ and $W$ is a map with the property
($\ref{eq:1}$). It results the Weyl connection $\tilde{\nabla}$
given by (\ref{eq:6}) without $Q$, which is a symmetric, compatible
linear connection. Supposing that the given distribution $\mathcal{D}$
is semi-Riemannian with respect to $g|_{\mathcal{D}}$ then it has
an orthogonal complementary distribution $\mathcal{D}^{\bot}$ with
the corresponding projector $Q^{\bot}$. Therefore we get two Weyl
substructures $(g_{|\mathcal{D}},W)$ and $(g_{|\mathcal{D}{}^{\bot}},W)$
with corresponding Weyl connections $\nabla$ and $\nabla^{\bot}$.
Using the terminology of \cite[p. 96]{b:f} let call $\nabla$ \textit{the
intrinsic Weyl connection of} $\mathcal{D}$ and $\nabla^{\bot}$
\textit{the transversal Weyl connection of} $\mathcal{D}.$ Using
the formula $(2.4)$ from \cite[p. 7]{b:f} it result a linear connection
$\nabla^{*}$ on $M$: \begin{equation}
\nabla_{X}^{*}Y=\nabla_{X}QY+\nabla_{X}^{\bot}Q^{\bot}Y,\label{eq:8}\end{equation}
 for $X,\, Y\in\chi(M).$ This linear connection is \textit{adapted}
to $\mathcal{D}$ and $\mathcal{D}{}^{\bot},$ i.e. for any $X\in\chi(M)$
and $U\in\Gamma(\mathcal{D})$,~$V\in\Gamma(\mathcal{D}{}^{\bot}),$
we have $\nabla_{X}^{*}U\in\Gamma(\mathcal{D})$, $\nabla_{X}^{*}V\in\Gamma(\mathcal{D^{\perp}})$
\cite[p. 7]{b:f}. The above formulae yield: \begin{equation}
\nabla_{X}^{*}Y=Q\tilde{\nabla}_{QX}QY+Q^{\bot}\tilde{\nabla}_{Q^{\perp}X}Q^{\bot}Y+Q[Q^{\bot}X,QY]+Q^{\bot}[QX,Q^{\bot}Y],\label{eq:9}\end{equation}
 which compared with relation $(3.16)$ from \cite[p. 17]{b:f} gives
a similar result to Theorem 5.3. of \cite[p. 26]{b:f} namely that
$\nabla^{*}$ is just \textit{the Vranceanu connection} defined by
the Weyl connection $\tilde{\nabla}$.
\end{example}
There are several features of the Vranceanu connection which makes
it important:

\begin{itemize}
\item it is defined on the sections of whole $TM$ not only of $\mathcal{D};$ 
\item if $\mathcal{D}$ is the tangent distribution of a foliation $\mathcal{F}$
(this case will be studied in the next section) then $\nabla^{*}$
is symmetric (torsion-free) if and only if the distribution $\mathcal{D}\mathcal{}^{\bot}$
is integrable, Theorem 1.5. of \cite[p. 100]{b:f};
\item if $\nabla^{*}$ is symmetric then the almost product structure $P=Q-Q^{\bot}$,
naturally associated to the decomposition (\ref{eq:3}) is integrable,
which means that the Nijenhuis tensor field $N_{P}$ vanishes.\\

\end{itemize}
Let us end this section with another form of the compatibility condition,
more precisely one in terms of quasi-connections \cite[p. 660]{s:h}.
Let $F\in\mathcal{T}_{1}^{1}(M)$ be a tensor field of $(1,1)$-type. 

\begin{defn}
An application $D:\,\chi(M)\times\Gamma(\mathcal{D})\rightarrow\Gamma(\mathcal{D})$
is a \emph{quasi-connection with respect to $F$} on $\mathcal{D}$
if, for all $X,Y\in\chi(M)$ and $Z\in\Gamma(\mathcal{D}):$

i) $D_{fX+gY}Z=fD_{X}Z+gD_{Y}Z,\,\, D_{X+Y}Z=D_{X}Z+D_{Y}Z,$

ii) $D_{X}(fZ)=fD_{X}Z+FX(f)Z.$\\

Let us remark that a linear connection $\nabla$ on $\mathcal{D}$
yields a quasi-connection $D^{\nabla}$ through \cite[p. 660]{s:h}:

\begin{equation}
D_{X}^{\nabla}Z=\nabla_{FX}Z\label{eq:10}\end{equation}
 and then we get:
\end{defn}
\begin{prop}
A linear connection is compatible with the Weyl substructure $W$
if and only if the associated quasi-connection (\ref{eq:10}) with
respect to the projector $Q$ makes $g$ a recurrent tensor with the
recurrence factor $-W(g)\circ Q$. 
\end{prop}
\medskip{}

\section{The compatible connection in foliated manifolds}

Let $(M,g)$ be an $(n+p)$-dimensional semi-Riemannian manifold and
$\mathcal{F}$ be an $n$-foliation on $M$. We assume that $\mathcal{D}$,
the tangent distribution of $\mathcal{F}$, is semi-Riemannian that
is, the induced metric tensor field on $\mathcal{D}$ is non-degenerate
and with constant index. The complementary orthogonal distribution
$\mathcal{D}^{\bot}$ to $\mathcal{D}$ is semi-Riemannian too, \cite[p. 95]{b:f};
let call $\mathcal{D}$ and $\mathcal{D}^{\bot}$ \textit{the structural}
and \textit{transversal distribution} respectively. Now, we want an
expression of the compatible connection in local coordinates.\\
So, let $\{\frac{\partial}{\partial x^{i}},\frac{\partial}{\partial x^{\alpha}}\}$
be a frame field adapted to the decomposition: \begin{equation}
TM=\mathcal{D}\oplus\mathcal{D}^{\bot},\label{eq:11}\end{equation}
 i.e. $i\in\{1,...,n\},\alpha\in\{n+1,...,n+p\}$ and $\frac{\partial}{\partial x^{i}}\in\Gamma(\mathcal{D}).$
With: \begin{equation}
g_{ij}=g\left(\frac{\partial}{\partial x^{i}},\frac{\partial}{\partial x^{j}}\right),g_{i\alpha}=g\left(\frac{\partial}{\partial x^{i}},\frac{\partial}{\partial x^{\alpha}}\right)\label{eq:12}\end{equation}
 it results an adapted basis for $\mathcal{D}^{\bot}$, \cite[p. 98]{b:f}:
\begin{equation}
\frac{\delta}{\delta x^{\alpha}}=\frac{\partial}{\partial x^{\alpha}}-A_{\alpha}^{i}\frac{\partial}{\partial x^{i}},\label{eq:13}\end{equation}
 where: \begin{equation}
A_{\alpha}^{i}=g^{ij}g_{j\alpha}.\label{eq:14}\end{equation}
 Remark that $\{\frac{\delta}{\delta x^{\alpha}}\}$ is orthogonal
to $\{\frac{\partial}{\partial x^{i}}\}.$ Let us point that $A_{\alpha}^{i}g_{i\beta}=A_{\beta}^{j}g_{j\alpha}.$
With respect to this adapted frame field we set: \begin{equation}
\nabla_{\frac{\partial}{\partial x^{j}}}\frac{\partial}{\partial x^{i}}=C_{ij}^{k}\frac{\partial}{\partial x^{k}},\nabla_{\frac{\delta}{\delta x^{\alpha}}}\frac{\partial}{\partial x^{i}}=D_{i\alpha}^{k}\frac{\partial}{\partial x^{k}}.\label{eq:15}\end{equation}

Then a computation similar to that of \cite[p. 99-100]{b:f} yields:

\medskip{}

\begin{prop}
The local coefficients of the compatible connection $\nabla$ with
respect to the semi-holonomic frame field $\{\frac{\partial}{\partial x^{i}},\frac{\delta}{\delta x^{\alpha}}\}$
are given by: \begin{eqnarray}
C_{ij}^{k} & =\Gamma_{ij}^{k}+\frac{1}{2}(\theta_{i}\delta_{j}^{k}+\theta_{j}\delta_{i}^{k}-g_{ij}\theta^{k}), & \ \ D_{i\alpha}^{k}=\frac{\partial A_{\alpha}^{k}}{\partial x^{i}}\label{eq:16}\end{eqnarray}
 where: \begin{equation}
W(g)=\theta_{i}\delta x^{i}+\rho_{\alpha}dx^{\alpha},\label{eq:17}\end{equation}
 $A_{\alpha}^{i}$ are given by (\ref{eq:14}), $\theta^{k}=g^{kl}\theta_{l}$,
$\Gamma_{ij}^{k}$ are the Christoffel symbols of $g$ with respect
to $\mathcal{D}$: \begin{equation}
\Gamma_{ij}^{k}=\frac{1}{2}g^{kl}\left(\frac{\partial g_{lj}}{\partial x^{i}}+\frac{\partial g_{il}}{\partial x^{j}}-\frac{\partial g_{ij}}{\partial x^{l}}\right)\label{eq:18}\end{equation}
 and $\{\delta x^{i},dx^{\alpha}\}$ is the dual frame of the given
semi-holonomic frame, with: \begin{equation}
\delta x^{i}=dx^{i}+A_{\alpha}^{i}dx^{\alpha}.\label{eq:19}\end{equation}

\end{prop}
For $p=0$ the well-known expression of \cite{f:l} and \cite{Sen}
are reobtained.

\medskip{}

\begin{example}
Let us continue Example 1.4 with $\mathcal{D}$ the tangent distribution
of a foliation $\mathcal{F}.$ Consider the local expressions above
to which we add: \begin{equation}
g_{\alpha\beta}=g\left(\frac{\delta}{\delta x^{\alpha}},\frac{\delta}{\delta x^{\beta}}\right)\label{eq:20}\end{equation}
 and denotes by $[g^{\lambda\mu}]$ the inverse matrix of $[g_{\alpha\beta}]$.
We express the Vranceanu connection $\nabla^{*}$ in local coordinates:
\begin{equation}
\left\{ \begin{array}{ll}
\nabla_{\frac{\partial}{\partial x^{j}}}^{*}\frac{\partial}{\partial x^{i}}= & C_{ij}^{k}\frac{\partial}{\partial x^{k}},\quad\nabla_{\frac{\delta}{\delta x^{\alpha}}}^{*}\frac{\partial}{\partial x^{i}}=D_{i\alpha}^{k}\frac{\partial}{x^{k}}\\
\\\nabla_{\frac{\partial}{\partial x^{i}}}^{*}\frac{\delta}{\delta x^{\alpha}}= & \nabla_{\frac{\partial}{\partial x^{i}}}^{\bot}\frac{\delta}{\delta x^{\alpha}}=L_{\alpha i}^{\gamma}\frac{\delta}{\delta x^{\gamma}},\ \\
\\\nabla_{\frac{\delta}{\delta x^{\beta}}}^{*}\frac{\delta}{\delta x^{\alpha}}= & \nabla_{\frac{\delta}{\delta x^{\beta}}}^{\bot}\frac{\delta}{\delta x^{\alpha}}=F_{\alpha\beta}^{\gamma}\frac{\delta}{\delta x^{\gamma}}.\end{array}\right.\label{eq:21}\end{equation}

\end{example}
A similar calculus like in \cite[p. 99]{b:f} gives \[
L_{\alpha i}^{\gamma}=0\]
 and: \begin{eqnarray}
F_{\alpha\beta}^{\gamma} & = & \frac{1}{2}g^{\gamma\mu}\left(\frac{\delta g_{\mu\beta}}{\delta x^{\alpha}}+\frac{\delta g_{\alpha\mu}}{\delta x^{\beta}}-\frac{\delta g_{\alpha\beta}}{\delta x^{\mu}}\right)\label{eq:22}\\
 &  & +\frac{1}{2}\left(\rho_{\alpha}\delta_{\beta}^{\gamma}+\rho_{\beta}\delta_{\alpha}^{\gamma}-\rho^{\gamma}g_{\alpha\beta}\right),\nonumber \end{eqnarray}
 with $\rho^{\gamma}=\rho_{\alpha}g^{\alpha\gamma}$.

\medskip{}

From Proposition 1.4. of \cite[p. 100]{b:f} it results that the only
non-null component of the torsion tensor field $T^{*}$ of $\nabla^{*}$
is: \begin{equation}
T^{*}\left(\frac{\delta}{\delta x^{\beta}},\frac{\delta}{\delta x^{\alpha}}\right)=T_{\alpha\beta}^{*k}\frac{\partial}{\partial x^{k}}=\left[\frac{\delta}{\delta x^{\alpha}},\frac{\delta}{\delta x^{\beta}}\right]=\left(\frac{\delta A_{\alpha}^{k}}{\delta x^{\beta}}-\frac{\delta A_{\beta}^{k}}{\delta x^{\alpha}}\right)\frac{\partial}{\partial x^{k}},\label{eq:23}\end{equation}
 which, therefore, describes exactly how far is $\mathcal{D}^{\perp}$
from integrability.

\medskip{}

Now, we consider the curvature tensor field $R^{*}$ of $\nabla^{*}$
and take in attention the transversal part using the notation of \cite[p. 104]{b:f08}:
\begin{equation}
\left\{ \begin{array}{llll}
R^{*}\left(\frac{\delta}{\delta x^{\gamma}},\frac{\delta}{\delta x^{\beta}}\right)\frac{\delta}{\delta x^{\alpha}} & = & R_{\alpha\beta\gamma}^{*\mu}\frac{\delta}{\delta x^{\mu}}\\
\\R^{*}\left(\frac{\partial}{\partial x^{i}},\frac{\delta}{\delta x^{\beta}}\right)\frac{\delta}{\delta x^{\alpha}} & = & R_{\alpha\beta i}^{*\mu}\frac{\delta}{\delta x^{\mu}}\\
\\R^{*}\left(\frac{\partial}{\partial x^{j}},\frac{\partial}{\partial x^{i}}\right)\frac{\delta}{\delta x^{\alpha}} & = & R_{\alpha ij}^{*\mu}\frac{\delta}{\delta x^{\mu}}\end{array}\right.\label{eq:24}\end{equation}
 with: \begin{equation}
\left\{ \begin{array}{ll}
R_{\alpha\beta\gamma}^{*\mu} & =\frac{\delta F_{\alpha\beta}^{\mu}}{\delta x^{\gamma}}-\frac{\delta F_{\alpha\gamma}^{\mu}}{\delta x^{\beta}}+F_{\alpha\beta}^{\varepsilon}F_{\varepsilon\gamma}^{\mu}-F_{\alpha\gamma}^{\varepsilon}F_{\varepsilon\beta}^{\mu}\\
\\R_{\alpha\beta i}^{*\mu} & =\frac{\partial F_{\alpha\beta}^{\mu}}{\partial x^{i}},\quad R_{\alpha ij}^{*\mu}=0.\end{array}\right.\label{eq:25}\end{equation}

The structural components of the curvature of $\nabla^{*}$ are: \begin{equation}
\left\{ \begin{array}{llll}
R^{*}(\frac{\delta}{\delta x^{\alpha}},\frac{\delta}{\delta x^{\beta}})\frac{\partial}{\partial x^{i}} & = & R_{i\alpha\beta}^{*h}\frac{\partial}{\partial x^{h}}\\
\\R^{*}(\frac{\partial}{\partial x^{k}},\frac{\delta}{\delta x^{\alpha}})\frac{\partial}{\partial x^{i}} & = & R_{i\alpha k}^{*h}\frac{\partial}{\partial x^{h}}\\
\\R^{*}(\frac{\partial}{\partial x^{k}},\frac{\partial}{\partial x^{j}})\frac{\partial}{\partial x^{i}} & = & R_{ijk}^{*h}\frac{\partial}{\partial x^{h}}\end{array}\right.\label{eq:26}\end{equation}
 with, \cite[p. 100]{b:f}: \begin{equation}
\left\{ \begin{array}{llll}
R_{i\alpha\beta}^{*h} & = & \frac{\delta D_{i\beta}^{h}}{\delta x^{\alpha}}-\frac{\delta D_{i\alpha}^{h}}{\delta x^{\beta}}+D_{i\beta}^{k}D_{k\alpha}^{h}-D_{i\alpha}^{k}D_{k\beta}^{h}-T_{\alpha\beta}^{*k}C_{ik}^{h}\\
\\R_{i\alpha k}^{*h} & = & \frac{\partial D_{i\alpha}^{h}}{\partial x^{k}}-\frac{\delta C_{ik}^{h}}{\delta x^{\alpha}}+D_{i\alpha}^{j}C_{jk}^{h}-C_{ik}^{j}D_{j\alpha}^{h}+D_{k\alpha}^{j}C_{ij}^{h}\\
\\R_{ijk}^{*h} & = & \frac{\partial C_{ij}^{h}}{\partial x^{k}}-\frac{\partial C_{ik}^{h}}{\partial x^{j}}+C_{ij}^{l}C_{lk}^{h}-C_{ik}^{l}C_{lj}^{h}.\end{array}\right.\label{eq:27}\end{equation}

Let $X\in\chi(M)$ with the decomposition $X=X^{i}\frac{\partial}{\partial x^{i}}+X^{\alpha}\frac{\delta}{\delta x^{\alpha}}$.
The covariant derivative of the metric $g$ with respect to the Vranceanu
connection is: \begin{equation}
\left\{ \begin{array}{llll}
(\nabla_{X}^{*}g)(\frac{\partial}{\partial x^{i}},\frac{\partial}{\partial x^{j}}) & = & X^{k}\left(\frac{\partial g_{ij}}{\partial x^{k}}-C_{ik}^{h}g_{hj}-C_{kj}^{h}g_{ih}\right)+X^{\alpha}\frac{\delta g_{ij}}{\delta x^{\alpha}}\\
\\(\nabla_{X}^{*}g)(\frac{\delta}{\delta x^{\alpha}},\frac{\delta}{\delta x^{\beta}}) & = & X^{i}\frac{\partial g_{\alpha\beta}}{\partial x^{i}}+X^{\mu}\left(\frac{\delta g_{\alpha\beta}}{\delta x^{\mu}}-F_{\alpha\mu}^{\rho}g_{\rho\beta}-F_{\mu\beta}^{\rho}g_{\alpha\rho}\right)\\
\\(\nabla_{X}^{*}g)(\frac{\partial}{\partial x^{i}},\frac{\delta}{\delta x^{\alpha}}) & = & 0\end{array}\right.\label{eq:28}\end{equation}
 and a straightforward computations using (\ref{eq:22}) yields: \begin{equation}
(\nabla_{X}^{*}g)\left(\frac{\delta}{\delta x^{\alpha}},\frac{\delta}{\delta x^{\beta}}\right)=X^{i}\frac{\partial g_{\alpha\beta}}{\partial x^{i}}-(X^{\mu}\rho_{\mu})g_{\alpha\beta}\label{eq:29}\end{equation}
 which implies a generalization of equivalence of items (i) and (ii)
of Theorem 3.3. from \cite[p. 112]{b:f} (obtained for $W(g)=0$):

\medskip{}

\begin{prop}
Let $(M,g,\mathcal{F})$ be a semi-Riemannian manifold, where $\mathcal{F}$
is a non-degenerate foliation. Then $g$ is a bundle-like metric for
$\mathcal{F}$ if and only if there exists an 1-form $W(g)$ on $M$
such that the induced metric $g$ on $\mathcal{D}^{\perp}$ is a recurrent
tensor with respect to the Vranceanu connection of the Weyl manifold
$(M,g,W:g\rightarrow W(g)),$ with the recurrence factor $-W(g)\circ Q^{\perp}:$
\[
(\nabla_{X}^{*}g)(Q^{\perp}Y,Q^{\perp}Z)=-W(g)(X)g(Q^{\perp}Y,Q^{\perp}Z),\,\forall X,Y,Z\in\chi(M).\]

\end{prop}
\medskip{}

Also: \begin{equation}
(\nabla_{X}^{*}g)\left(\frac{\partial}{\partial x^{i}},\frac{\partial}{\partial x^{j}}\right)=-(X^{k}\theta_{k})g_{ij}+X^{\alpha}\frac{\delta g_{ij}}{\delta x^{\alpha}}.\label{eq:30}\end{equation}

\section{Weyl structures on tangent bundles of Finsler manifolds}

Let $N$ be a real $n$-dimensional manifold and $TN$ the tangent
bundle of $N$. Then a local chart $x=(x^{a})$ on $N$ defines a
local chart $(x,y)=(x^{a},y^{a})_{1\leq a\leq n}$ on $TN$. Denote
by $0(N)$ the zero section of $TM$ and consider $TN^{0}=TN\setminus0(N)$.

\medskip{}

\begin{defn}
The pair $(N,\, F)$ is a \textit{Finsler manifold} if \\
 $F:TN\rightarrow[0,\infty)$ with the following conditions: \\
 F1) $F$ is smooth on $TN^{0}$ and vanishes only on $0(N)$, \\
 F2) $F$ is positively homogeneous of degree one with respect to
$(y^{a})$, \\
 F3) the matrix $[g_{bc}(x^{a},y^{a})]=[\frac{1}{2}\frac{\partial^{2}F^{2}}{\partial y^{b}\partial y^{c}}]$
is positive definite.
\end{defn}
The vertical bundle $V(N)$ of $N$ is the tangent distribution to
the foliation defined by the fibers of $\pi:TN\rightarrow N.$ Then
$V(N)$ is locally spanned by ${\frac{\partial}{\partial y^{a}}}$.
Denote by $[g^{bc}]$ the inverse matrix of $[g_{bc}]$ and define:
\begin{equation}
G^{a}(x,y)=\frac{1}{4}g^{ab}\left(\frac{\partial^{2}F^{2}}{\partial y^{b}\partial x^{c}}y^{c}-\frac{\partial F^{2}}{\partial x^{b}}\right).\label{eq:31}\end{equation}

There exists on $TN$ an $n$-distribution $H(N)$, called \textit{horizontal},
locally spanned by the vector fields: \begin{equation}
\frac{\delta}{\delta x^{a}}=\frac{\partial}{\partial x^{a}}-G_{a}^{b}\frac{\partial}{\partial y^{b}},\label{eq:32}\end{equation}
 where: \begin{equation}
G_{b}^{a}=\frac{\partial G^{a}}{\partial y^{b}}.\label{eq:33}\end{equation}

It is easy to see that $H(N)$ is complementary to $V(N)$ in $TN$
and using the decomposition: \begin{equation}
T(TN)=H(N)\oplus V(N)\label{eq:34}\end{equation}
 we define the Riemannian metric $G$ on $TN$, called the \textit{Sasaki-Finsler
metric}: \begin{equation}
G=\left(\begin{array}{ll}
g_{ab} & 0\\
0 & g_{ab}\end{array}\right)\label{eq:35}\end{equation}
 which means that with respect to the semi-holonomic frame field $\{\frac{\delta}{\delta x^{a}},\frac{\partial}{\partial y^{a}}\}$
we have: \begin{equation}
G\left(\frac{\delta}{\delta x^{a}},\frac{\delta}{\delta x^{b}}\right)=G\left(\frac{\partial}{\partial y^{a}},\frac{\partial}{\partial y^{b}}\right)=g_{ab},G\left(\frac{\delta}{\delta x^{a}},\frac{\partial}{\partial y^{b}}\right)=0.\label{eq:36}\end{equation}

The above discussion shows that on the Riemannian manifold $(TN,G)$
we have a foliation $\mathcal{F}$ with $V(N)$ and $H(N)$ as structural
and transversal distribution respectively, therefore we obtain the
framework discussed in the previous section. Suppose given an 1-form
$W(G)$ on $TN$ with the expression \begin{equation}
W(G)=\rho_{a}dx^{a}+\theta_{a}\delta y^{a},\label{eq:37}\end{equation}
 where $(dx^{a},\delta x^{a})$ is the dual frame of the given semi-holonomic
frame: \begin{equation}
\delta y^{a}=dy^{a}+G_{b}^{a}dx^{b}.\label{eq:38}\end{equation}

The aim of this section is to obtain the coefficients of the Vranceanu
connection for the pair (Weyl structure $W:G\rightarrow W(G),$ distribution
$V(N)$). These coefficients are given by: \begin{equation}
\left\{ \begin{array}{ll}
\nabla_{\frac{\partial}{\partial y^{b}}}^{*}\frac{\partial}{\partial y^{a}}=C_{ab}^{c}\frac{\partial}{\partial y^{c}}, & \nabla_{\frac{\delta}{\delta x^{b}}}^{*}\frac{\partial}{\partial y^{a}}=D_{ab}^{c}\frac{\partial}{\partial y^{c}},\\
\\\nabla_{\frac{\partial}{\partial y^{b}}}^{*}\frac{\delta}{\delta x^{a}}=L_{ab}^{c}\frac{\delta}{\delta x^{c}}, & \nabla_{\frac{\partial}{\partial y^{b}}}^{*}\frac{\delta}{\delta x^{a}}=F_{ab}^{c}\frac{\delta}{\delta x^{c}}.\end{array}\right.\label{eq:39}\end{equation}

Using Proposition 3.1. from \cite[p. 226-227]{b:f} and the results
of the previous section we derive:

\medskip{}

\begin{prop}
The Vranceanu connection of a Weyl manifold $(TN,G,W)$ has the local
coefficients: \begin{equation}
\left\{ \begin{array}{llll}
C_{ab}^{c} & = & \frac{1}{2}\left(g^{cd}\frac{\partial g_{ab}}{\partial y^{d}}+\theta_{a}\delta_{b}^{c}+\theta_{b}\delta_{a}^{c}-\theta^{c}g_{ab}\right),\\
\\D_{ab}^{c} & = & \frac{\partial^{2}G^{a}}{\partial y^{b}\partial y^{c}},\,\,\, L_{ab}^{c}=0,\\
\\F_{ab}^{c} & = & \frac{1}{2}g^{cd}\left(\frac{\delta g_{db}}{\delta x^{a}}+\frac{\delta g_{ad}}{\delta x^{b}}-\frac{\delta g_{ab}}{\delta x^{d}}\right)+\frac{1}{2}\left(\rho_{a}\delta_{b}^{c}+\rho_{b}\delta_{a}^{c}-\rho^{c}g_{ab}\right).\end{array}\right.\label{eq:40}\end{equation}

\end{prop}
\medskip{}

\begin{example}
1) A Finsler manifold is a \textit{Landsberg space} if, \cite[p. 239]{a:t}:
\begin{equation}
\frac{1}{2}g^{cd}\left(\frac{\delta g_{db}}{\delta x^{a}}+\frac{\delta g_{ad}}{\delta x^{b}}-\frac{\delta g_{ab}}{\delta x^{d}}\right)=\frac{\partial G_{a}^{c}}{\partial y^{b}}.\label{eq:41}\end{equation}

\end{example}
Therefore, the Vranceanu connection for a Weyl manifold provided by
a Landsberg space has: \begin{equation}
F_{ab}^{c}=\frac{\partial^{2}G^{c}}{\partial y^{a}\partial y^{b}}+\frac{1}{2}\left(\rho_{a}\delta_{b}^{c}+\rho_{b}\delta_{a}^{c}-\rho^{c}g_{ab}\right).\label{eq:42}\end{equation}

2) A Finsler manifold is a \textit{locally Minkowski space} if, \cite[p. 239]{a:t},
there exists a covering by charts $(U,x)$ of $N$ such that $g_{ab}=g_{ab}(y)$.
A locally Minkowski space is a Landsberg one with $G^{a}=0$ and then
one have: \begin{equation}
F_{ab}^{c}=\frac{1}{2}\left(\rho_{a}\delta_{b}^{c}+\rho_{b}\delta_{a}^{c}-\rho^{c}g_{ab}\right).\label{eq:43}\end{equation}

\medskip{}

\begin{example}
In the following we consider a natural 1-form $W(G)$. The condition
F3 of Definition 3.1 means that $F^{2}:TN\rightarrow[0,\infty)$ is
a regular Lagrangian in the sense of Analytical Mechanics and then,
it defines a Legendre transform $L(F^{2}):TN\rightarrow T^{*}N$,
with $T^{*}N$ the cotangent bundle of $N$. With coordinates $(x^{a})$
on $TN$ we have induced coordinates $(x^{a},p_{a})$ on $T^{*}N$.
Also, on $T^{*}N$ lives a global 1-form, called Liouville, $\theta=p_{a}dx^{a}$.
The pullback of the Liouville form through the Legendre transform,
$\theta_{F}=L(F^{2})^{*}(\theta),$ is called \textit{the Cartan form}
of $F$. Therefore, we define $W(G)=\theta_{F}=\frac{1}{2}\frac{\partial F^{2}}{\partial y^{a}}dx^{a}$
which yields:
\end{example}
\medskip{}

\begin{prop}
The Vranceanu connection of a Weyl manifold \\
 $(TN,F,\theta_{F})$ has the coefficients: \begin{equation}
\left\{ \begin{array}{lll}
C_{ab}^{c} & = & \frac{1}{2}g^{cd}\frac{\partial g_{ab}}{\partial y^{d}},\\
\\D_{ab}^{c} & = & \frac{\partial^{2}G^{c}}{\partial y^{a}\partial y^{b}},\\
\\L_{ab}^{c} & = & 0,\\
\\4F_{ab}^{c} & = & 2g^{cd}\left(\frac{\delta g_{db}}{\delta x^{a}}+\frac{\delta g_{ad}}{\delta x^{b}}-\frac{\delta g_{ab}}{\delta x^{d}}\right)\\
\\ &  & +\frac{\partial F^{2}}{\partial y^{a}}\delta_{b}^{c}+\frac{\partial F^{2}}{\partial y^{b}}\delta_{a}^{c}-\frac{\partial F^{2}}{\partial y^{u}}g^{uc}g_{ab}.\end{array}\right.\label{eq:43}\end{equation}

\end{prop}
But $\frac{\partial F^{2}}{\partial y^{v}}=2g_{vu}y^{u}$ from F2
and then: \begin{equation}
F_{ab}^{c}=\frac{1}{2}g^{cd}\left(\frac{\delta g_{db}}{\delta x^{a}}+\frac{\delta g_{ad}}{\delta x^{b}}-\frac{\delta g_{ab}}{\delta x^{d}}\right)+\frac{1}{2}y^{u}\left(g_{ub}\delta_{a}^{c}+g_{au}\delta_{b}^{c}-g_{ab}\delta_{u}^{c}\right).\label{eq:44}\end{equation}

\begin{cor}
The Vranceanu connection of a Landsberg-Weyl manifold $(TN,F,\theta_{F})$
has: \begin{equation}
F_{ab}^{c}=\frac{\partial^{2}G^{c}}{\partial y^{a}\partial y^{b}}+\frac{1}{2}y^{u}\left(g_{ub}\delta_{a}^{c}+g_{au}\delta_{b}^{c}-g_{ab}\delta_{u}^{c}\right)\label{eq:45}\end{equation}
 while for a locally Minkowski space: \begin{equation}
F_{ab}^{c}=\frac{1}{2}y^{u}\left(g_{ub}\delta_{a}^{c}+g_{au}\delta_{b}^{c}-g_{ab}\delta_{u}^{c}\right).\label{eq:46}\end{equation}

\end{cor}
\begin{rem*}
From now we use $i,j,k,...$ for the vertical coordinates and $a,b,c,...$
for the horizontal ones, for a better identification of the structural
and respectevely transversal components. But we keep the above notations
for the coefficients of the Vranceanu connection.
\end{rem*}
The transversal components of the curvature $R^{*}$ of $\nabla^{*}$
are given by (\ref{eq:25}) with: \begin{equation}
R_{abi}^{*c}=\frac{1}{2}(g_{ib}\delta_{a}^{c}+g_{ai}\delta_{b}^{c}-g_{ab}\delta_{i}^{c})\label{eq:47}\end{equation}
 which never vanishes. 

\medskip{}

\begin{example}
Let $g=(g_{ab}(x))$ be a Riemannian metric on $N$ with the Christoffel
coefficients $\Gamma_{ab}^{c}$. Then $F=\left(g_{uv}y^{u}y^{v}\right)^{\frac{1}{2}}$
is a Finsler fundamental function on $N$.
\end{example}
\medskip{}

\begin{defn}
We define the \textit{Vranceanu-Cartan connection} on $TN$ for the
Riemannian manifold $(N,g)$, the Vranceanu connection obtained from
the process of Example 3.4. Namely it is associated to the Weyl manifold
$(TN,W:G\rightarrow\theta_{F})$ with the above $F$.
\end{defn}
\medskip{}

This Vranceanu connection is a particular case of Proposition 3.2
and then:

\medskip{}

\begin{prop}
The Vranceanu-Cartan connection of the tangent bundle $TN$ of a Riemannian
manifold $(N,g)$ is: \begin{equation}
\left\{ \begin{array}{ll}
C_{ab}^{c} & =0,\quad D_{ab}^{c}=\Gamma_{ab}^{c},\quad L_{ab}^{c}=0,\\
\\F_{ab}^{c} & =\Gamma_{ab}^{c}+\frac{1}{2}y^{u}\left(g_{ub}\delta_{a}^{c}+g_{au}\delta_{b}^{c}-g_{ab}\delta_{u}^{c}\right).\end{array}\right.\label{eq:48}\end{equation}
 For $X=X^{i}\frac{\partial}{\partial y^{i}}+X^{a}\frac{\delta}{\delta x^{a}}$
the non-null covariant derivatives of the Sasaki-Riemann metric $G$
with respect to the Vranceanu-Cartan connection are: \begin{equation}
\left\{ \begin{array}{ll}
(\nabla_{X}^{*}G)(\frac{\partial}{\partial y^{i}},\frac{\partial}{\partial y^{j}}) & =X^{a}\frac{\delta g_{ij}}{\delta x^{a}}=X^{a}\frac{\partial g_{ij}}{\partial x^{a}},\\
\\(\nabla_{X}^{*}G)(\frac{\delta}{\delta x^{a}},\frac{\delta}{\delta x^{b}}) & =-(g_{ci}X^{c}y^{i})g_{ab}.\end{array}\right.\label{eq:49}\end{equation}
\\

\end{prop}
\begin{rem}
Denoting by $Rg$ the $(1,3)$-Riemannian curvature tensor field of
$g$ and using (\ref{eq:25}), we get that the non-vanishing transversal
components of the curvature of Vranceanu - Cartan connection are: 
\end{rem}
\begin{itemize}
\item $R_{abc}^{*d}=(Rg)_{abc}^{d}+$ a very complicated expression in $g$
and $y$, \\

\item $R_{abi}^{*c}$ from (\ref{eq:47}), 
\end{itemize}
while the only non-null structural component of $R^{*}$ is, using
(\ref{eq:27}):

\begin{itemize}
\item $R_{iab}^{*j}=(Rg)_{iab}^{j}$. \\

\end{itemize}
Let us point also, that $T_{ab}^{*c}$ from (\ref{eq:23}) is, \cite[p. 233]{b:f}:

\[
T_{ab}^{*c}=(Rg)_{dab}^{c}y^{d}.\]
\\

Denote with $V$ and $H$ the vertical and horizontal projectors of
$TN.$ They correspond to $Q$ respectively $Q^{\bot}$ in the notations
of the first two sections. The above discussion about the curvature
of $\nabla^{*}$ yields:

\medskip{}

\begin{prop}
For the Vranceanu-Cartan connection and $X,Y,Z\in\chi(TN)$ we have:
\\
 1) $\nabla^{*}$ is torsion-free if and only if the base manifold
$(N,g)$ is flat.\\
 2) $R^{*}(HX,HY)VZ=0$ if and only if the base manifold $(N,g)$
is flat.\\
 Moreover, the Vranceanu-Cartan is never vertical-horizontal flat
but is vertical flat i.e $R^{*}(V\cdot,V\cdot)=0$.
\end{prop}
\medskip{}

Using the equivalent conditions from \cite[p. 237]{a:t} we derive:

\medskip{}

\begin{cor}
The projection $\pi_{T}:(TN,G)\rightarrow(N,g)$ is totally geodesic
i.e. the projection of any geodesic in $(TN,G)$ is also a geodesic
in $(N,g)$ if and only if the Vranceanu-Cartan connection is torsion-free.
\end{cor}
\medskip{}

Let us end this section with the covariant derivative of the Liouville
vector fields with respect to the Vranceanu connection in the general
(i.e. Finslerian) framework of this section. More precisely, define
\cite[p. 231]{b:f}: \begin{equation}
\left\{ \begin{array}{ll}
L & =y^{i}\frac{\partial}{\partial y^{i}},\\
\\L^{*} & =y^{a}\frac{\delta}{\delta y^{a}}\end{array}\right.\label{eq:51}\end{equation}
 called \textit{the Liouville vector field} on $TN$ respectively
the \textit{transversal Liouville vector field} or \textit{geodesic
spray}. For $X=X^{i}\frac{\partial}{\partial y^{i}}+X^{a}\frac{\delta}{\delta x^{a}}$
it results: \begin{equation}
\left\{ \begin{array}{ll}
\nabla_{X}^{*}L & =(X^{i}+X^{j}y^{k}C_{kj}^{i})\frac{\partial}{\partial y^{i}}+X^{c}(D_{bc}^{a}y^{b}-G_{c}^{a})\delta_{a}^{i}\frac{\partial}{\partial y^{i}},\\
\\\nabla_{X}^{*}L^{*} & =(X^{i}+X^{j}y^{k}L_{kj}^{i})\delta_{i}^{a}\frac{\delta}{\delta x^{a}}+X^{c}(F_{bc}^{a}y^{b}-G_{c}^{a})\frac{\delta}{\delta x^{a}}\end{array}\right.\label{eq:52}\end{equation}
 which, replacing the coefficients from (\ref{eq:40}) and using the
relations $(3.32_{a})$ from \cite[p. 231]{b:f} and $(3.39)$ from
\cite[p. 232]{b:f}, yields:

\medskip{}

\begin{prop}
The covariant derivative of the Liouville vector fields with respect
to the Vranceanu connection of a Weyl manifold $(TN,G,W)$ are: \begin{equation}
\left\{ \begin{array}{ll}
\nabla_{X}^{*}L & =\left[X^{i}+\frac{1}{2}X^{j}y^{k}(\theta_{j}\delta_{k}^{i}+\theta_{k}\delta_{j}^{i}-\theta^{i}g_{jk})\right]\frac{\partial}{\partial y^{i}},\\
\\\nabla_{X}^{*}L^{*} & =\left[X^{i}\delta_{i}^{a}+\frac{1}{2}X^{b}y^{c}\left(\rho_{b}\delta_{c}^{a}+\rho_{c}\delta_{b}^{a}-\rho^{a}g_{bc}\right)\right]\frac{\delta}{\delta x^{a}}.\end{array}\right.\label{eq:53}\end{equation}
 In particular, for the Weyl manifold $(TN,F,\theta_{F})$ we get:
\begin{equation}
\left\{ \begin{array}{ll}
\nabla_{X}^{*}L & =X^{i}\frac{\partial}{\partial y^{i}},\\
\\\nabla_{X}^{*}L^{*} & =\left(X^{i}\delta_{i}^{a}+\frac{1}{2}X^{a}y^{c}y^{u}g_{uc}\right)\frac{\delta}{\delta x^{a}}.\end{array}\right.\label{eq:54}\end{equation}

\end{prop}
In order to provide a global expression of these relations let us
denote the vertical and horizontal components of $W(G)$: \begin{equation}
W(G)^{V}=\theta_{i}\delta y^{i},\quad W(G)^{H}=\rho_{a}dx^{a}\label{eq:55}\end{equation}
 and then (\ref{eq:53}) becomes: \[
\left\{ \begin{array}{ll}
2\nabla_{X}^{*}L & =2VX+W(G)^{V}(VX)L+W(G)^{V}(L)VX-G(VX,L)W(G)^{V\#}\\
\\2\nabla_{X}^{*}L^{*} & =2\Theta(X)+W(G)^{H}(HX)L^{*}+W(G)^{H}(L^{*})HX-G(HX,L^{*})W(G)^{H\#}\end{array}\right.\]
 where $\#$ is the musical isomorphism defined by $G$: \begin{equation}
W(G)^{V\#}=\theta^{i}\frac{\partial}{\partial y^{i}},\quad W(G)^{H\#}=\rho^{a}\frac{\delta}{\delta x^{a}}\label{eq:56}\end{equation}
 while (\ref{eq:54}) is: \begin{equation}
\left\{ \begin{array}{ll}
\nabla_{X}^{*}L & =VX,\\
\\\nabla_{X}^{*}L^{*} & =\Theta(X)+\frac{1}{2}F^{2}HX\end{array}\right.\label{eq:57}\end{equation}
 with $F$ the Finsler fundamental function of Definition 3.1 and
$\Theta$ the \textit{adjoint structure,} \cite[p. 991]{m:ab}, $\Theta=\frac{\delta}{\delta x^{i}}\otimes\delta y^{i}$.

\medskip{}

\vspace{0.3cm}

\curraddr{\noindent Oana Constantinescu \quad{}Mircea Crasmareanu \\
 Faculty of Mathematics \\
 University \char`\"{}Al. I. Cuza\char`\"{} \\
 Ia\c{s}i, 700506,\\
 Romania }

\email{\noindent E-mail: oanacon@uaic.ro \quad{}mcrasm@uaic.ro}

\medskip{}

\email{\noindent http://www.uaic.ro/$\sim$mcrasm\quad{}http://www.uaic.ro/$\sim$oanacon}
\end{document}